\documentclass[12pt,draft]{article}

\usepackage{amsfonts}
\usepackage{amsmath}
\usepackage{amssymb}
\usepackage{amscd}
\usepackage{amsthm}
\usepackage{indentfirst}
\usepackage{epsfig}

\usepackage[vmargin=4.4cm,hmargin=3.6cm]{geometry}

\newtheorem{thm}{Theorem}[section]
\newtheorem{lemma}[thm]{Lemma}
\newtheorem{prop}[thm]{Proposition}
\newtheorem{cor}[thm]{Corollary}

\theoremstyle{definition}

\newtheorem{defin}[thm]{Definition}
\newtheorem{rem}[thm]{Remark}

\newtheorem*{prf}{Proof}
\newtheorem*{notation}{Notation}
\newtheorem*{acknow}{Acknowledgements}

\newcommand{\R}{{\mathbb{R}}}

\newcommand{\Z}{{\mathbb{Z}}}

\newcommand{\cI}{{\mathcal{I}}}

\newcommand{\cM}{{\mathcal{M}}}

\newcommand{\cO}{{\mathcal{O}}}

\newcommand{\cT}{{\mathcal{T}}}
\newcommand{\cU}{{\mathcal{U}}}

\newcommand{\cC}{{\mathcal{C}}}

\newcommand{\fc}{{:\ }}
\newcommand{\ve}{\varepsilon}

\newcommand{\tb}{\textbf}

\newcommand{\ol}{\overline}

\DeclareMathOperator{\id}{id}

\DeclareMathOperator{\Hom}{Hom}

\DeclareMathOperator{\End}{End}

\DeclareMathOperator{\Crit}{Crit}

\DeclareMathOperator{\ind}{ind}
\DeclareMathOperator{\Int}{Int}
\DeclareMathOperator{\codim}{codim}

\DeclareMathOperator{\im}{im}

\begin{document}


\title{Functoriality in Morse theory on closed manifolds}
\author{Avraham Aizenbud and Frol Zapolsky}
\maketitle

\begin{abstract}
We develop functoriality for Morse theory, namely, to a pair of Morse-Smale systems and a generic smooth map between the underlying manifolds we associate a chain map between the corresponding Morse complexes, which descends to the correct map on homology. This association does not in general respect composition. We give sufficient conditions under which composition is preserved. As an application we provide a new proof that the cup product as defined in Morse theory on the chain level agrees with the cup product in singular cohomology. In appendices we present a proof (due to Paul Biran) that the unstable manifolds of a Morse-Smale system are the open cells of a CW structure on the underlying manifold, and also we show that the Morse complex of the triple is canonically isomorphic to the cellular complex of the CW structure. This gives a new proof that the Morse complex is actually a complex and that it computes the homology of the manifold.
\end{abstract}

\section{Introduction}

Morse theory attaches to a generic triple \((M,F,\rho)\), where \(M\) is a closed manifold, \(F\) a Morse function on \(M\) and \(\rho\) is a Riemannian metric on \(M\), a chain complex \(\cM_*(M,F,\rho)\), called the Morse complex of the triple (also known as the Thom-Smale-Witten complex in the literature). The homology of this complex is known to be canonically identified with the singular homology of \(M\), and so is independent of \((F,\rho)\), see Appendix \ref{app_Morse_complex}. For general background on Morse homology, see for example \cite{Morse_theory_Austin}, \cite{Morse_homology_Schwarz}. The Morse complex itself does depend on this data. In this paper we investigate some functorial properties of this dependence. Namely, given two such triples, we associate to a generic smooth map between the underlying manifolds a chain map between the Morse complexes, which descends to the correct map on homology. However, this association is in general not functorial, that is, it does not respect composition. We show that under some condition on the maps in question composition is preserved. Also, we believe that in general some weaker form of functoriality is present. At this point, we do not know what is the appropriate framework for this question.

\subsection{Structure of the paper}

In Section \ref{Definitions} we fix notations, give main definitions and formulate the results. Section \ref{construct_proofs} is dedicated to proofs of the results. In addition we apply functoriality to prove that the cup product on Morse complexes descends to the usual cup product on cohomology. This definition of the cup product in Morse homology has been around for some time now, but there seems to be no quotable reference for the details. In any case, we show how it fits into the general theory in a more coherent fashion.

In Appendix \ref{app_transversality} we prove that transverse maps (see definition \ref{map_types}) form an open and dense subset in the set of smooth maps. Transverse maps are the smooth maps to which we associate a chain map between the Morse complexes. In Appendix \ref{app_intersect_orient} we treat the relevant part of intersection theory of (oriented) manifolds. In Appendix \ref{app_CW_struct} we bring an expanded version of the proof by Paul Biran (see \cite{Biran_lagrangian_barriers_symp_embeddings}) of the fact that the unstable manifolds of a Morse-Smale triple are the open cells of a CW structure on the underlying manifold and in Appendix \ref{app_Morse_complex} we show that the Morse complex of the triple is canonically isomorphic to the cellular complex of the CW complex. Note that this gives an alternative proof of the fact that the Morse complex is indeed a chain complex and that its homology computes the homology of the manifold. In a sense, this proof is more intuitive than the ones which are usually encountered. Its drawback is its heavy reliance on difficult theorems of the theory of dynamical systems.

\begin{acknow} We would like to thank Paul Biran, Michael Farber, Dmitry Gourevitch, Michael Hutchings, Matthias Kreck, Janko Latschev, Dimitri Novikov, Leonid Polterovich, Peter Pushkar, and Eugenii Shustin for useful discussions.
\end{acknow}

\section{Definitions}\label{Definitions}

\subsection{Preliminaries on Morse theory}

Let \(M\) be a closed manifold, \(F \in C^\infty(M)\) be a Morse function, and \(\rho\) a Riemannian metric on \(M\). If \(x\) is a critical point of \(F\), then the Hessian \(d^2_xF\) is a well-defined nondegenerate quadratic form on \(T_xM\) and the dimension of the negative space of \(d^2_xF\) is called the index of \(F\) at \(x\), denoted by \(\ind_F x\). The set of critical points of \(F\) of index \(k\) is \(\Crit_k(F)\), while the set of all critical points \(\Crit(F) = \bigcup_{k=0}^n \Crit_k(F)\), where \(n = \dim M\).

The negative gradient flow of the triple \((M,F,\rho)\) is the flow of the vector field \(-\nabla_\rho F\). Morse functions will be denoted by variations of the letter \(F\), while the corresponding flows will be denoted by the suitable variation of the small letter \(f\) with superscript \(t\), for example \(F_2 \rightsquigarrow f_2^t\).

Let \((M,F,\rho)\) be a triple as above, and let \(x \in \Crit_k(F)\). The \emph{unstable manifold} through \(x\) is
\[U^x := \big\{y \in M\,|\, \lim_{t \to - \infty}f^t(y) = x\big\}\,,\]
and the \emph{stable manifold} through \(x\) is
\[S^x := \big\{y \in M\,|\, \lim_{t \to \infty}f^t(y) = x\big\}\,.\]
It is a standard fact that \(U^x\) and \(S^x\) are diffeomorphic to \(\R^k\) and \(\R^{n-k}\), respectively. The triple \((M,F,\rho)\) is called \emph{Morse-Smale} if every stable manifold intersects every unstable manifold transversely.

\begin{defin}
Let \((M,F,\rho)\) be a Morse-Smale triple. There is a natural CW structure on \(M\) corresponding to the triple, the open cells being the unstable manifolds, see Appendix \ref{app_CW_struct}. Let \(M^k\), for \(k = 0,\dots,\dim M\), denote the \(k\)-th skeleton of this structure, that is, \(M^k\) is the union of the unstable manifolds of dimension up to \(k\).
\end{defin}

\begin{rem}\label{CW_char_maps}
A CW complex for us is a CW complex with fixed attaching maps. In particular, each open cell is endowed with a smooth structure and an orientation. However, in case the CW complex in question comes from a Morse-Smale triple, the smooth structure on an open cell is that of the unstable manifold, not the one induced from the attaching map, which only provides an orientation for the cell. The characteristic maps, the smooth structures, and the orientations obtained by this convention are henceforth referred to as canonical.
\end{rem}

\subsection{Definitions and results}

\begin{notation} We shall use the following notational convention. Whenever there is given a map of sets \(\varphi \fc M \to N\) and a family of self-maps of \(N\), \(\{f^t\}_{t \in \R}\), we shall denote \(\varphi^t := f^t \circ \varphi\). Now if \((N,F,\rho)\) is a Morse-Smale triple, there is a canonical family of self-maps of \(N\), namely the negative gradient flow \(f^t\), and it is this family that will always be understood in the notation \(\varphi^t\).
\end{notation}

\begin{defin}\label{map_types}
Let \(M_1\) be a CW complex and let \((M_2,F_2,\rho_2)\) be a Morse-Smale triple. Let \(\varphi \fc M_1 \to M_2\) be a continuous map. Define the map \(\varphi^\infty(x) := \lim_{t\to\infty}\varphi^t(x)\) (\(\varphi^\infty\) is in general discontinuous). The map \(\varphi\) is called:
\begin{itemize}
\item\emph{regular}, if \(\varphi^\infty(M_1^k) \subset M_2^k\) for all \(k\);

\item\emph{cellwise smooth}, or in short, \emph{cw-smooth}, if the restriction of \(\varphi\) to every open cell of \(M_1\) (taken with the canonical smooth structure, see Remark \ref{CW_char_maps}) is smooth;

\item\emph{transverse}, if it is cw-smooth, regular, and if for any open cell \(U_1 \subset M_1\) of dimension \(k\) and any stable manifold \(S_2 \subset M_2\) of codimension \(k\) the restriction \(\varphi|_{U_1}\) is transverse to \(S_2\).
\end{itemize}
\end{defin}

\begin{rem}
If \((M_1,F_1,\rho_1)\) is a Morse-Smale triple and \(N\) is a manifold, then a smooth map \(\varphi \fc M_1 \to N\) is of course cw-smooth for the canonical CW structure on \(M_1\), see Remark \ref{CW_char_maps}.
\end{rem}

\begin{rem}
In the definition of a transverse map the triple \((M_2,F_2,\rho_2)\) needs not be Morse-Smale. A triple \((M,F,\rho)\) is Morse-Smale if and only if the identity \(\id_M\) is a transverse map.
\end{rem}


\begin{defin}For a Morse-Smale triple $(M,F,\rho)$ we let $\cM_*(M,F,\rho)$ be the Morse complex of the triple, $\cM_k(M,F,\rho)$ being generated by the critical points of $F$ of index $k$.

\end{defin}

\begin{defin}
Let \(M_1\) be a CW complex and let \((M_2,F_2,\rho_2)\) be a Morse-Smale triple. Let \(\varphi \fc M_1 \to M_2\) be a transverse map. For any \(k\) let \(U_1 \subset M_1\) be a \(k\)-cell and let \(S_2 \equiv S^{x_2} \subset M_2\) be the stable manifold corresponding to a critical point \(x_2 \in \Crit_k(F_2)\). Then we can define the number
\[n_\varphi(U_1,x_2) := \#\varphi|_{U_1} \cap S_2\,,\]
that is the intersection index. For the precise definition, as well as various orientation issues, we refer the reader to Appendix \ref{app_intersect_orient}. Using these numbers, we can define a sequence of homomorphisms (for notations concerning CW complexes, see Subsection \ref{metric_CW})
\[\cM(\varphi)_k \fc \cC_k(M_1) \to \cM_k(M_2,F_2,\rho_2)\,,\]
as follows. For an open cell \(U_1 \in \cC_k(M_1)\) put
\[\cM(\varphi)_k(U_1) := \sum_{x_2 \in \Crit_k(F_2)} n_\varphi(U_1,x_2)x_2\,.\]
\end{defin}

Now we are ready to formulate the main results.

\begin{thm}\label{transverse_funct}
Let \(M_1\) be a CW complex and let \((M_2,F_2,\rho_2)\) be a Morse-Smale triple. Let \(\varphi \fc M_1 \to M_2\) be a transverse map. Then \(\cM(\varphi)\) is a chain map, and the induced map \(\cM(\varphi)_*\) on homology coincides with \(\varphi_*\), the homology of $\cM_*(M_2,F_2,\rho_2)$ being identified with the singular homology of $M_2$.
\end{thm}

\begin{lemma}\label{transversality}
Let \((M_i,F_i,\rho_i)\), \(i=1,2\) be Morse-Smale triples. The set of transverse maps \(M_1 \to M_2\) is open and dense in \(C^\infty(M_1,M_2)\).
\end{lemma}

The proof is given in Appendix \ref{app_transversality}.

\begin{thm} \label{C0_extension}
The association \(\varphi \mapsto \cM(\varphi)\) has a unique \(C^0\)-stable extension to the set of regular maps between a CW complex and a Morse-Smale triple. This extension still associates a chain map to a regular map, and for any regular \(\varphi\) we have \(\cM(\varphi)_* = \varphi_*\).
\end{thm}

Theorems \ref{transverse_funct} and \ref{C0_extension} are proved in Section \ref{construct_proofs}.

This association is not in general functorial, for two reasons. Firstly, the composition of two regular maps needs not be regular, and secondly, even when \(\cM\) is defined, it does not in general commute with composition. There is however, a certain amount of functoriality possessed by \(\cM\), and the following definition and theorem cover some of it.

\begin{defin}
Let \((M_1,F_1,\rho_1)\) and \((M_2,F_2,\rho_2)\) be Morse-Smale triples. A regular map \(\psi \fc M_1 \to M_2\) is called \emph{stably regular} if for any \(x \in M_1\) there is \(y \in \Crit(F_2)\) such that \(\psi(f_1^t(x)) \in S^y\) for all \(t \in \R\).
\end{defin}

\begin{thm}\label{stable_regular_funct}
Let \(M_1\) be a CW complex, and let \((M_i,F_i,\rho_i)\), \(i=2,3\) be Morse-Smale triples. Let \(\varphi \fc M_1 \to M_2\) be a regular map, and let \(\psi \fc M_2 \to M_3\) be a stably regular map. Then \(\psi \circ \varphi\) is regular, and moreover
\[\cM(\psi \circ \varphi) = \cM(\psi) \circ \cM(\varphi)\,.\]
\end{thm}

The proof is given in Subsection \ref{stably_regular_maps}.

\section{Constructions and proofs}\label{construct_proofs}

\subsection{Generalized Lyapunov theorem}

The following proposition says that if a closed subset is carried by the negative gradient flow pointwise close to a skeleton, then it also happens uniformly.

\begin{prop}\label{generalized_Lyapunov}Let \((M,F,\rho)\) be a Morse-Smale triple, and let \(X\) be a subcomplex of \(M\).

(i) Let \(V \supset X\) be an open neighborhood. Then there is an open neighborhood \(U \supset X\) such that \(U \subset V\) and \(U\) is invariant under the flow \(f^t\), that is \(f^t(U) \subset U\) for all \(t \geq 0\);

(ii) Let \(Y \subset M\) be a closed subset. Suppose that for any \(y \in Y\) and any open \(V \supset X\) there exists \(T > 0\) such that \(f^T(y)\in V\). Then for any open \(V \supset X\) there is \(T > 0\) such that \(f^T(y) \in V\) for any \(y \in Y\).
\end{prop}

\begin{prf}We prove both statements simultaneously by cellular induction. The case of points is clear. Let us prove (i). Let \(X'\) be a subcomplex of \(M\), and let \(X = X' \cup e\), where \(e\) is a cell not in \(X'\), such that \(\partial e \subset X'\) and \(e\) is the unstable manifold of a critical point \(x\). Let \(U' \supset X'\) be an \(\{f^t\}_{t\geq0}\)-invariant open neighborhood such that \(U' \subset V\). Let \(\Phi \fc (-1,1)^n \to M\) be a coordinate chart with \(\Phi(0)=x\) and \(\Phi(x_1,\dots,x_k,0,\dots,0)\in U^x\) for \(x_i\in(-1,1)\). For \(r \in (0, 1)\) let
\[C=\{(x_1,...,x_n) \in (-1,1)^n | \sum_{i=1}^k x_i^2=1/4,\, \sum_{i=k+1}
^n x_i^2 = 0\}\]

\[S(r)=\{(x_1,...,x_n) \in (-1,1)^n | \sum_{i=1}^k x_i^2=1/4,\, \sum_{i=k+1}
^n x_i^2 \leq r^2\}\]
and
\[B(r)=\{(x_1,...,x_n) \in (-1,1)^n | \sum_{i=1}^k x_i^2\leq 1/4,\, \sum_{i=k+1}
^n x_i^2 \leq r^2\}\,.\]
By the inductive assumption there is \(T>0\) such that \(f^T(\Phi(C)) \subset U'\). Hence there is \(r>0\) such that \(f^T(\Phi(S(r)))\subset U'\). Let \(\Psi \fc S(r) \times [0, T] \to M\) be the flow map. By the compactness of \([0, T]\) there is \(r' > 0\) such that \(\Psi(S(r') \times [0, T]) \subset V\). Now define \(U = \Int\big(\Phi(B(r')) \cup \Psi(S(r') \times [0, T])\big) \cup U'\). Then \(U\) satisfies the requirements of (i).

(ii) By (i) there is an open \(U \subset V\) which is invariant under the flow. Let \(y \in Y\). There is \(T_y > 0\) such that \(f^{T_y}(y) \in U\). Since \(U\) is open, there is an open neighborhood \(W_y\) of \(y\) such that \(f^{T_y}(W_y) \subset U\). The sets \(W_y\) form an open cover of \(Y\), and owing to the compactness of \(Y\) there is a finite subcover \(\{W_{y_i}\}_i\). Now let \(T = \max_i T_{y_i}\). \qed
\end{prf}

\subsection{Metric CW complexes}\label{metric_CW}

Let \((X,d)\) be a (pseudo-)metric space of finite diameter. If \(Z\) is a set, then on the set of maps \(\{f\fc Z \to X\}\) there is a natural (pseudo-)metric, which we shall denote by \(d\) as well, defined by \(d(f,g):= \sup_{z \in Z}d(f(z),g(z))\). For \(\delta > 0\) and \(A \subset X\) we denote \(A_\delta = \{x \in X \,|\, d(x,A) < \delta\}\). Given \(\ve > 0\), a homotopy \(\{\varphi_t \fc Y \to X\}_{t \in J}\), where \(J \subset \R\) is some interval and \(Y\) is a topological space, is called an \(\ve\)-\emph{homotopy}, if for all \(t,t' \in J\) we have \(d(\varphi_t,\varphi_{t'}) < \ve\).

In this subsection all CW complexes are finite, in particular, compact and metrizable.

\begin{defin}
A pseudo-metric CW complex is a pair \((M,d)\) where \(M\) is a CW complex and \(d\) is a pseudo-metric which induces a topology coarser than that of \(M\). The pair \((M,d)\) is a metric CW complex if \(d\) is a metric inducing the topology of \(M\).
\end{defin}

Recall that to a CW complex \(M\) one associates its chain complex \(\cC(M)\) and to any cellular map \(\varphi \fc M \to N\) a chain map \(\cC(\varphi) \fc \cC(M) \to \cC(N)\). If \(\varphi_t\) is a homotopy consisting of cellular maps, then \(\cC(\varphi_t)\) is independent of \(t\).

\begin{defin}
Let \(M\) be a CW complex and let \((N,d)\) be a metric CW complex. A continuous map \(\varphi \fc M \to N\) is called \(\ve\)-cellular, where \(\ve > 0\), if for any integer \(k \geq 0\) we have \(\varphi(M^k)  \subset (N^k)_\ve\).
\end{defin}

It is well-known that any map between CW complexes is homotopic to a cellular map. We shall make use of a refinement of this statement:

\begin{prop}Let \(N\) be a CW complex and let \((M,d)\) be a metric CW complex. Then for any \(\ve > 0\) there exists \(\delta > 0\) such that for any \(\delta\)-cellular map \(\varphi \fc N \to M\) there is an \(\ve\)-homotopy \(\{\varphi_t\}_{t \in I}\) with \(\varphi_0 = \varphi\) and \(\varphi_1\) cellular. \qed
\end{prop}

This follows from quantitative cellular approximation, see Theorem \ref{quant_cell_approx} below. The proof of the next result goes along the same lines, and will be omitted.

\begin{prop}Let \(M\) be a CW complex and let \((N,d)\) be a metric CW complex. Then for any \(\ve > 0\) there exists \(\delta > 0\) such that for any homotopy \(\varphi_t \fc M \to N\) of \(\delta\)-cellular maps there exists a homotopy \(\psi_t\) of cellular maps such that \(d(\varphi_t,\psi_t)<\ve\) for any \(t\). \qed
\end{prop}

The following definitions are used in the next subsection.

\begin{defin}
Let \((M,d)\) be a metric CW complex. To it we associate a number \(\eta_1 \equiv \eta_1(M,d) > 0\) such that for any CW complex \(L\) and any pair of cellular maps \(\varphi_1,\varphi_2 \fc L \to M\) with \(d(\varphi_1,\varphi_2) < \eta_1\) we have \(\cC(\varphi_1) = \cC(\varphi_2)\). The existence of \(\eta_1\) is evident.
\end{defin}

\begin{defin}
Let \((M,d)\) be a metric CW complex. To it we associate a number \(\eta_2 \equiv \eta_2(M,d) > 0\) such that for any CW complex \(L\) and any \(\eta_2\)-cellular map \(\varphi \fc L \to M\) there exists a cellular map \(\psi\) with \(d(\varphi,\psi)< \frac {\eta_1}2\).
\end{defin}

\begin{defin}
Let \(N\) be a CW complex and let \((M,d)\) be a metric CW complex. Put \(\eta_2 := \eta_2(M,d)\). Let \(\varphi \fc N \to M\) be \(\eta_2\)-cellular. Choose a cellular map \(\psi\) with \(d(\varphi,\psi)< \frac {\eta_1(M,d)}2\). Define \(\cC(\varphi) := \cC(\psi)\). Clearly this is independent of \(\psi\). If \(\varphi\) is a cellular map, then this definition coincides with the original one.
\end{defin}

The rest of this subsection is dedicated to the proof of Theorem \ref{quant_cell_approx}. The proof is along the lines of that of the classical cellular approximation theorem, see for example \cite{alg_topology_Hatcher}.

\begin{lemma} Let \((N,d)\) be a metric CW complex and let \(K \subset N\) be a subcomplex. Consider the pseudo-metric CW complex \((N \times I,d')\), where  \(d'\) is given by \(d'\big((z,t),(z',t')\big):=d(z,z')\). Then for any \(\ve > 0\) there exists a retract \(r \fc N \times I \to N \times 0 \cup K \times I \subset N \times I\) with \(d'(r,\id_{N \times I}) < \ve\). \qed
\end{lemma}

\begin{lemma}\label{extending_homotopies} Let \((M,D)\) be a metric space. Let \(N\) be a CW complex, and let \(K \subset N\) be a subcomplex. Suppose that \(p \fc N \to M\) is continuous, and \(\{q_t \fc K \to M\}_{t \in I}\) is an \(\ve\)-homotopy with \(q_0 = p|_K\). Then there exists a \(2\ve\)-homotopy \(\{p_t \fc N \to M\}\) such that \(p_0 = p\) and \(p_t|_K = q_t\).
\end{lemma}

\begin{prf} Metrize \(N\) by a metric \(d\). Since \(N\) is compact and \(p\) is continuous, it is uniformly continuous, that is, there is \(\delta > 0\) such that \(d(z,z') < \delta \Rightarrow D(p(z),p(z')) < \ve\). By the previous lemma there is a retract \(r \fc N \times I \to N\times 0 \cup K \times I\) such that \(d(r,\id_{N \times I}) < \delta\). Put \(p_t(z) = \overline p(r(z,t))\), where \(\overline p \fc N \times 0 \cup K \times I \to M\) is the union of the maps \(p\) and \(q_t\).
\qed
\end{prf}

\begin{lemma}\label{quant_cell_approx_one_cell} Let \((M,d)\) be a metric CW complex. Then for any \(\ve > 0\) there exists \(\delta>0\) such that if \(N\) is a CW complex, \(U \subset N\) is a top-dimensional cell, and \(\varphi \fc N \to M\) is a \(\delta\)-cellular map, which is cellular on \(N - U\), then there is an \(\ve\)-homotopy \(\varphi_t\), stationary on \(N-U\), with \(\varphi_0 = \varphi\) and \(\varphi_1\) cellular.
\end{lemma}

For the proof we need the following straightforward

\begin{lemma} Let \(D \subset \R^n\) be the standard closed ball of radius \(1\). Let \(d\) be a pseudo-metric on \(D\), which is continuous with respect to the standard metric, and whose restriction to the interior of \(D\) induces the standard topology. Then for any \(\ve > 0\) there is \(\delta > 0\) such that there exists an \(\ve\)-homotopy \(r_t \fc D \to D\), such that \(r_t|_{\partial D} = \id\), \(r_0 = \id_D\) and \(r_1((\partial D)_\delta) = \partial D\) (the \(\delta\)-neighborhood is taken with respect to \(d\)). \qed
\end{lemma}

\begin{prf}[of Lemma \ref{quant_cell_approx_one_cell}] Let \(k = \dim U\). Proceed by induction on the maximal dimension \(l\) (which we of course suppose to be \(> k\)) of a cell of \(M\) which meets \(\varphi(U)\). Let \(s\) be the number of cells in \(M^l-M^k\) which meet \(U\). Let \(V \subset M\) be such a cell, and let \(\Phi \fc D \to \overline V\) be the characteristic map. Let \(d\) be the pseudo-metric on \(D\) induced by \(\Phi\). For \(\ve'=\frac{\ve}{2s}\) choose \(\delta > 0\) and an \(\ve'\)-homotopy \(r_t\) as in the previous lemma. Define \(\overline r_t \fc M^l \to M^l\) by \(\overline r_t|_{M^l - V} = \id\), \(\overline r_t(\Phi(z)) := \Phi(r_t(z))\) for \(z \in D\). Repeat the process for each cell of dimension \(l\). Concatenate all the obtained homotopies into one \(\ve/2\)-homotopy, spanned over a suitable interval of time, denote it by \(\overline r_t^l\). Compose: \(\overline r_t^l \circ \varphi\). The end point \(\varphi'\) of this homotopy sends \(U\) into \(M^{l-1}\), and we can apply the inductive assumption to get a further \(\ve/2\)-homotopy whose endpoint sends \(\varphi'(U)\) into \(M^k\). The concatenation of these two \(\ve/2\)-homotopies is the desired \(\ve\)-homotopy. \qed
\end{prf}

\begin{thm}[Quantitative cellular approximation theorem]\label{quant_cell_approx} Let \((M,d)\) be a metric CW complex, and let \(N\) be a CW complex. For any \(\ve > 0\) there exists \(\delta > 0\) such that if \(K \subset N\) is a subcomplex and \(\varphi \fc N \to M\) is a \(\delta\)-cellular map, which is cellular on \(K\), then there is an \(\ve\)-homotopy \(\varphi_t\), stationary on \(K\), with \(\varphi_0 = \varphi\) and \(\varphi_1\) cellular.
\end{thm}

\begin{prf}
By cellular induction on \(N\). Let \(U \subset N - K\) be a cell. By the inductive assumption, for any \(\ve_1 > 0\) there exists \(\delta > 0\) such that for any \(\delta\)-cellular map \(N - U \to M\) which is cellular on \(K\) there is an \(\ve_1\)-homotopy, stationary on \(K\), between it and a cellular map.

Let \(\delta_1 > 0\) be such that for any \(\delta_1\)-cellular map \(\varphi \fc N \to M\) which is cellular on \(N-U\), there exists an \(\ve/2\)-homotopy \(\varphi_t\) fixed on \(N-U\) such that \(\varphi_0 = \varphi\), and \(\varphi_1\) is cellular.

Let \(\ve_1 = \min(\delta_1/4,\ve/4)\) and let \(\delta < \delta_1/2\) be as provided by the inductive assumption (see the first paragraph), corresponding to \(\ve_1\).

Let \(\varphi \fc N \to M\) be a \(\delta\)-cellular map which is cellular on \(K\), and let \(\varphi' = \varphi|_{N-U}\). Let \(\varphi''_t\) be an \(\ve_1\)-homotopy between \(\varphi'\) and a cellular map \(N-U \to M\), provided by the inductive assumption. Extend \(\varphi''_t\) to a \(2\ve_1\)-homotopy \(\varphi'_t \fc N \to M\). Let \(\psi = \varphi'_1\). Since \(\varphi'_0 = \varphi\) is \(\delta\)-cellular and \(\varphi'_t\) is a \(2\ve_1\)-homotopy, \(\psi\) is \((2\ve_1 + \delta)\)-cellular, in particular, it is \(\delta_1\)-cellular, since \(2\ve_1+\delta < \delta_1\).

Since \(\psi\) is \(\delta_1\)-cellular and \(\psi|_{N-U}\) is cellular, there exists an \(\ve/2\)-homotopy \(\psi_t\) on \(N\), stationary on \(N-U\), such that \(\psi_0 = \psi\) and \(\psi_1\) is cellular. The concatenation of the \(\ve/2\)-homotopies \(\varphi'_t\) (\(\varphi'_t\) is a \(2\ve_1\)-homotopy, hence an \(\ve/2\)-homotopy) and \(\psi_t\) is the desired \(\ve\)-homotopy. \qed
\end{prf}

\subsection{Regular maps}

In this subsection a Morse-Smale triple \((M,F,\rho)\) is viewed as a metric CW complex, where the CW structure is the canonical one (see Remark \ref{CW_char_maps}), and the metric, denoted by \(d\), is the one induced by \(\rho\). For the purposes of this subsection we need to modify the definition of the number \(\eta_2(M)\) in case \(M\) is the underlying manifold of a Morse-Smale triple.

\begin{defin}
Let \((M,F,\rho)\) be a Morse-Smale triple. Let \(\eta_2:=\eta_2(M,d)\) (\(d\) being the Riemannian distance). Denote by \(\eta_2'\equiv \eta_2'(M,F,\rho)\) some number \(\leq \eta_2\) such that for any \(k < l\) and any \(y \in \Crit_l(F)\) the stable manidold \(S^y\) does not meet the \(\eta_2'\)-neighborhood of \(M^k\).
\end{defin}

\begin{defin}
Let \(M_1\) be a CW complex and let \((M_2,F_2,\rho_2)\) be a Morse-Smale triple. A map \(\varphi \fc M_1 \to M_2\) is called \emph{\(t\)-admissible}, where \(t > 0\), if \(\varphi^{s}\) is \(\eta_2'(M_2,F_2,\rho_2)\)-cellular for any \(s \geq t\).
\end{defin}

\begin{rem}
A map is regular if and only if it is \(t\)-admissible for some \(t\), due to Lemma \ref{generalized_Lyapunov}.
\end{rem}

\begin{defin}
Let \(M_1\) be a CW complex and let \((M_2,F_2,\rho_2)\) be a Morse-Smale triple. For a \(t\)-admissible \(\varphi \fc M_1 \to M_2\) we define \(\widetilde \cM_s(\varphi):= \cC(\varphi^s)\) for \(s \geq t\).
\end{defin}

The association \([t,\infty) \ni s \mapsto \widetilde \cM_s (\varphi)\) is clearly locally constant, hence constant, and we have

\begin{lemma}
If \(\varphi\) is regular, then \(\widetilde \cM_s(\varphi)\) is independent of \(s\) for \(s \geq t\), where \(t\) is such that \(\varphi\) is \(t\)-admissible. We denote by \(\widetilde \cM (\varphi)\) the common value. \qed
\end{lemma}

\begin{defin}
Let \(M_1\) be a CW complex and let \((M_2,F_2,\rho_2)\) be a Morse-Smale triple. A map \(\varphi \fc M_1 \to M_2\) is called \emph{\(\ve\)-stable} if for any map \(\varphi'\) with \(d(\varphi,\varphi')<\ve\) we have: (i) \(\varphi'\) is regular and (ii) \(\widetilde\cM(\varphi') = \widetilde\cM(\varphi)\).
\end{defin}

\begin{lemma}\label{t_admissible_ve_stable}
Let \(M_1\) be a CW complex and let \((M_2,F_2,\rho_2)\) be a Morse-Smale triple. For any \(t > 0\) there is \(\ve > 0\) such that any \(t\)-admissible map \(\varphi \fc M_1 \to M_2\) is \(\ve\)-stable.
\end{lemma}

\begin{prf}
This is due to the uniform continuity of the map \(f^t\). \qed
\end{prf}

\begin{cor}
\(\widetilde\cM\) is \(C^0\)-stable. \qed
\end{cor}

\subsection{\(\widetilde\cM = \cM\)}

In this subsection we prove Theorems \ref{transverse_funct} and \ref{C0_extension}. They are both corollaries of the following

\begin{thm}\label{trans_reg_compatibility}
Let \(M_1\) be a CW complex and let \((M_2,F_2,\rho_2)\) be a Morse-Smale triple. Let \(\varphi \fc M_1 \to M_2\) be a transverse map. Then \(\cM(\varphi) = \widetilde\cM(\varphi)\).
\end{thm}

\begin{prf}
If \(S\) is a metric sphere, there exists a positive number, which we denote \(\eta(S)\), such that any two continuous maps from a topological space to \(S\) which are \(\eta(S)\)-close, are homotopic.

For any \(x \in \Crit_k(F_2)\) choose a parametrization \(\Phi_x \fc B_1 \times B_2 \to M_2\), where \(B_1 \subset \R^k,B_2\subset\R^{n-k}\) are small balls centered at \(0\), \(n = \dim M_2\), such that \(\Phi_x(0) = x\), $\Phi(B_1\times\{0\})\subset U^x$, $\Phi(\{0\}\times B_2)\subset S^x$. Denote $T_x = \Phi(B_1\times B_2)$, $D_x = \Phi(B_1\times \{0\})$ and let $\pi_x\fc T_x \to D_x$ be defined by $\pi_x(\Phi(a,b))=\Phi(a,0)$. Let \(S_x = M_2^k/(M_2^k - D_x)\).

Let
\[\ve = \frac 1 2 \min_{x \in \Crit(F_2)} \eta(S_x)\,.\]
Choose \(\delta > 0\) such that (i) any \(\delta\)-cellular map from \(M_1\) to \(M_2\) can be approximated by an \(\ve\)-close cellular map, (ii) \((M_2^k - D_x)_\delta \cup T_x \supset (M_2^k)_\delta\) for all \(k\) and \(x \in \Crit_k(F_2)\), (iii) \(\delta < \eta_2(M_2)\), (iv) for any \(y \in T_x \cap (M_2^k)_\delta\) we have \(d(\pi_x(y),y) < \ve\). Choose \(T > 0\) such that for any \(t \geq T\)  the map \(\varphi^t\) is \(\delta\)-cellular.

Let \(\psi \fc M_1 \to M_2\) be a cellular map such that \(d(\psi,\varphi^T) < \ve\).

Fix an open \(k\)-dimensional cell \(U_1 \subset M_1\) and \(x \in \Crit_k(F_2)\). We shall define \(\varphi' \fc M_1^k/(M_1^k-U_1) \to S_x\) as follows. The point \([M_1^k - U_1]\) is mapped to \([M_2^k-D_x]\). Now let \(y \in U_1\). If \(\varphi^T(y) \in T_x\), then \(\varphi'(y) = \pi_x(\varphi^T(y))\). If \(\varphi^T(y) \notin T_x\), \(\varphi'(x) = [M_2^k-D_x]\). The map \(\varphi'\) is clearly continuous. Let \(\psi' \fc M_1^k/(M_1^k-U_1) \to M_2^k/(M_2^k-U^x)\) be the map induced by \(\psi\) . Define \(\psi'' = p \circ \psi'\), where \(p \fc M_2^k/(M_2^k-U^x) \to S_x\) is the projection. Clearly \(d(\psi'',\varphi') < 2\ve\), hence \(\deg \psi'' = \deg \varphi'\). On the other hand, \(\deg \psi' = \deg \psi''\), since \(p\) is a homotopy equivalence. In addition, \(x \in S_x\) is a regular value of \(\varphi'\), hence \(\deg \varphi' = n_\varphi(U_1,x)\).

To summarize, \(\deg \psi' = n_\varphi(U_1,x)\), hence \(\cM(\varphi) = \widetilde \cM(\varphi)\). \qed

\end{prf}

Henceforth for a regular map \(\varphi\) we shall denote \(\cM(\varphi) := \widetilde\cM(\varphi)\), and thanks to the theorem, this notation is consistent.

\subsection{Stably regular maps}\label{stably_regular_maps}

In this subsection we prove Theorem \ref{stable_regular_funct}. Throughout let \(M_1\) be a CW complex and \(\{(M_i,F_i,\rho_i)\}_{i=2,3}\) be a pair of Morse-Smale triples.

\begin{prop}\label{homotopy_thru_regular} Let \(\varphi_t \fc M_1 \to M_2\) be a homotopy through regular maps. Then \(\cM(\varphi_t)\) is independent of \(t\). \qed
\end{prop}

\begin{lemma}\label{regular_comp_cellular}
If \(\varphi \fc M_1 \to M_2\) is cellular, \(\psi \fc M_2 \to M_3\) is regular, then \(\psi \circ \varphi\) is regular and \(\cM(\psi \circ \varphi) = \cM(\psi) \circ \cC(\varphi)\). 
\qed
\end{lemma}

\begin{prop}\label{composition_large_t}
Let \(\varphi \fc M_1 \to M_2\) and \(\psi \fc M_2 \to M_3\) be regular maps. Then for all large enough \(t\) the composition \(\psi \circ \varphi^t\) is regular and \(\cM(\psi \circ \varphi^t)=\cM(\psi) \circ \cM(\varphi)\).
\end{prop}

We need the following

\begin{lemma}
There exist \(t_0\) and \(\ve > 0\) such that for any \(t \geq t_0\) the composition \(\psi \circ \varphi^t\) is \(\ve\)-stable.
\end{lemma}

\begin{prf}
Thanks to Lemma \ref{t_admissible_ve_stable}, if we show that there are \(t_0\) and \(T\) such that \(\psi \circ \varphi^t\) is \(T\)-admissible for all \(t \geq t_0\), we shall obtain the desired assertion.

Choose \(T \geq 0\) such that for all \(k\) the set \(\psi^T(M_2^k)\) lies in an invariant neighborhood \(V^k\) of \(M_3^k\), contained in the \(\eta_2'(M_3)\)-neighborhood of \(M_3^k\). Let \(U^k\) be an invariant neighborhood of \(M_2^k\) such that \(\psi^T(U^k) \subset V^k\). Now choose \(t_0\) such that \(\varphi^{t_0}(M_1^k) \subset U^k\). \qed
\end{prf}

\begin{prf}[of the proposition]
Only the second assertion is nontrivial. Let \(t_0\) and \(\ve\) be as provided by the lemma. Since \(\psi\) is uniformly continuous, there is \(\delta > 0\) such that for any \(x,y\) with \(d(x,y) < \delta\) we have \(d(\psi(x),\psi(y))<\ve\). There is \(t_1 > t_0\) such that for any \(t > t_1\) there is a cellular map \(\varphi'\) such that \(d(\varphi^t,\varphi') < \delta\), and \(\cM(\varphi) = \cC(\varphi')\). It follows that \(d(\psi \circ \varphi^t, \psi \circ \varphi') < \ve\), hence \(\cM(\psi \circ \varphi^t) = \cM(\psi \circ \varphi')\). Now use Lemma \ref{regular_comp_cellular} to conclude that \(\cM(\psi \circ \varphi^t) = \cM(\psi \circ \varphi') = \cM(\psi)\circ\cC(\varphi') = \cM(\psi)\circ\cM(\varphi)\). \qed
\end{prf}

\begin{lemma} Let \(\varphi \fc M_1 \to M_2\) be regular and \(\psi \fc M_2 \to M_3\) be stably regular. Then \(\psi \circ \varphi\) is regular.
\end{lemma}

\begin{prf}
Let \(x \in M_1^k\). We need to show that \(\psi(\varphi(x))\) lies in some \(S^z\), where \(z \in \Crit_l(F_3)\) with \(l \leq k\). Since \(\psi\) is regular, \(\psi(M_2^k) \subset \bigcup_{z \in \Crit_l(F_3), l \leq k}S^z\). It follows that there is \(\delta > 0\) such that the whole \(\delta\)-neighborhood of \(M_2^k\) is sent \(\psi\) into the latter union. Since \(\varphi\) is regular, there is \(t\) such that \(d(\varphi^t(x),M_2^k) < \delta\), thus \(\psi(\varphi^t(x)) \in S^z\), where \(z \in \Crit(F_3)\) of index \(\leq k\). By the definition of stable regularity, this means that the whole open gradient line through \(\varphi^t(x)\), which is the same as the gradient line through \(\varphi(x)\), is sent by \(\psi\) into \(S^z\), that is \(\psi(\varphi(x)) \in S^z\). \qed
\end{prf}

\begin{prf}[of Theorem \ref{stable_regular_funct}]
We only need to show that \(\cM(\psi \circ \varphi) = \cM(\psi) \circ \cM(\varphi)\). Since \(\varphi\) is regular, so is \(\varphi^t\) for any \(t\), that is \(\{\psi \circ \varphi^t\}_t\) is a homotopy through regular maps, and by Proposition \ref{homotopy_thru_regular} we have \(\cM(\psi \circ \varphi) = \cM(\psi \circ \varphi^0) = \cM(\psi \circ \varphi^t)\). By Proposition \ref{composition_large_t} \(\cM(\psi \circ \varphi^t) = \cM(\psi) \circ \cM(\varphi)\) for all large \(t\). \qed
\end{prf}

\subsection{Cup product}

\begin{defin}
Let \((M,F,\rho)\) be a Morse-Smale triple. Denote
\[\cM^k (M,F,\rho) := \Hom(\cM_k (M,F,\rho),\Z)\,.\]
Together with the differential \(\partial f:= f \circ \partial\), where \(f \in \cM^*\) and \(\partial\) is the Morse differential on \(\cM_*\), this constitutes a cochain complex, whose cohomology is canonically identified with the singular cohomology of \(M\). For a CW complex \(M_1\), a Morse-Smale triple \((M_2,F_2,\rho_2)\), and a regular map \(\varphi \fc M_1 \to M_2\) define
\[\cM(\varphi)^k \fc \cM^k(M_2,F_2,\rho_2) \to \cC^k(M_1)\]
by \(\cM(\varphi)^k = (\cM(\varphi)_k)^*,\)
the star meaning the dual map. This is clearly a cochain map which, thanks to Theorem \ref{C0_extension}, gives the correct map on cohomology.
\end{defin}

We will now define the cup product on Morse cohomology using the functoriality, mimicking the analogous construction in the singular theory.

\begin{notation} Let \((M,F,\rho)\) be a Morse-Smale triple, and let \(x \in \Crit_k(F)\). Define \(\delta_x \in \cM^k(M,F,\rho)\) by \(\delta_x(x) = 1\) and \(\delta_x(y) = 0\) for \(y \neq x\).
\end{notation}

Let \(M\) be a closed manifold. Choose three Morse functions \(F_i\), \(i=1,2,3\) and three Riemannian metrics \(\rho_i\), \(i=1,2,3\), such that \((M,F_i,\rho_i)\) are Morse-Smale triples, and such that the diagonal map \(\Delta \fc (M,F_3,\rho_3) \to (M \times M,F_1+F_2,\rho_1+\rho_2)\) is transverse. Let \(x_i \in \Crit_{k_i}(F_i)\) be such that \(k_3 = k_1 + k_2\). Define the number \(n_{x_3}^{x_1,x_2} := \#U^{x_3} \cap S^{x_1} \cap S^{x_2}\), that is the triple intersection number, see Appendix \ref{app_intersect_orient}. Due to the fact that \(\Delta\) is a transverse map, this number is well-defined. Define a bilinear map
\[\smile \fc \cM^{k_1}(M,F_1,\rho_1) \times \cM^{k_2}(M,F_2,\rho_2) \to \cM^{k_3}(M,F_3,\rho_3)\]
by
\[(\delta_{x_1} \smile \delta_{x_2})(x_3) = n_{x_3}^{x_1,x_2}\,.\]
Theorem \ref{transverse_funct} yields that \(\smile\) is a bi-chain map, and the induced map on cohomology coincides with the cup product. This formula for the cup product has been known before.

\appendix

\section{Transversality}\label{app_transversality}

Here we prove Lemma \ref{transversality}. For this we need a couple of definitions and lemmas.

If \(M,N\) are closed manifolds, \(M' \subset M, N' \subset N\) are submanifolds, put
\[\cU(M',N') = \{\varphi \in C^\infty(M,N)\,|\, \varphi|_{M'} \text{ is transverse to } N'\}\,,\]
and
\[\cI(M',N') = \{\varphi \in C^\infty(M,N)\,|\, \varphi({\overline{M'}-M'})\cap N' = \varphi(M')\cap (\overline{N'}-N') = \varnothing\}\,.\]

\begin{lemma}The set \(\cU(M',N')\) is dense and the set \(\cU(M',N')\cap\cI(M',N')\) is open in \(C^\infty(M,N)\).
\end{lemma}

\begin{prf}
The first statement is a standard transversality argument, and will be omitted.

Let \(\varphi \in \cU(M',N')\cap\cI(M',N')\). It is clear that for \(x \in M',y\in N'\) there is an open neighborhood \(\cO_{x,y}\) of \(\varphi\) and open neighborhoods \(U_{x,y} \ni x,V_{x,y} \ni y\) such that for any \(\varphi' \in \cO\) we have that \(\varphi'|_{U_{x,y} \cap M'}\) is transverse to \(V_{x,y} \cap N'\). Also there is an open neighborhood \(\cO\) of \(\varphi\) and open neighborhoods \(U \supset \overline{M'} - M',V \supset \overline{N'} - N'\) such that for any \(\varphi' \in \cO\) we have \(\varphi'(U) \cap N' = \varphi'(M') \cap V = \varnothing\).

Now \(\{U_{x,y} \times V_{x,y}\}_{(x,y) \in M' \times N'} \cup \{U \times \overline{N'}, \overline{M'} \times V\}\) is an open cover of the compact set \(\overline{M'} \times \overline{N'}\), hence there are \(\{(x_i,y_i) \in M' \times N'\}_{i=1}^r\) such that \(\{U_{x_i,y_i} \times V_{x_i,y_i}\}_{i=1}^r \cup \{U \times \overline{N'}, \overline{M'} \times V\}\) is still a cover. Now \(\bigcap_{i=1}^r \cO_{x_i,y_i} \cap \cO\) is an open neighborhood of \(\varphi\) inside \(\cU(M',N')\cap\cI(M',N')\). \qed

\end{prf}

\begin{prf}[of Lemma \ref{transversality}]
The set of transverse maps equals
\[\cT = \bigcap \{\cU(U^x,S^y)\,|\,x\in\Crit_k(F_1),y\in\Crit_l(F_2),k\leq l\}\,.\]
Let \(\widetilde{M_2}\) be the CW complex whose cells are the stable manifolds of the triple \((M_2,F_2,\rho_2)\), or, which is the same, the CW complex of the triple \((M_2,-F_2,\rho_2)\).
For a subcomplex \(X \subset M_1 \times \widetilde{M_2}\) define
\[\cT(X) = \bigcap \{\cU(U^x,S^y)\,|\,(x,y)\in\Crit_k(F_1) \times \Crit_l(F_2) \cap X,k\leq l\}\,.\]
Clearly \(\cT = \cT(M_1\times \widetilde{M_2})\). We shall prove by cellular induction on \(X\) that \(\cT\) is open and dense. For \(X = \varnothing\) the statement is trivial. Assume the statement for a certain \(X\) and let \(X' = X \uplus e\), where \(e = U^x \times S^y\) is a cell not in \(X\). We have
\[\cT(X') = \cT(X) \cap \cU(U^x,S^y)\,.\]
Now \(\cT(X)\) is open dense by inductive assumption, \(\cU(U^x,S^y)\) is dense by the lemma, hence \(\cT(X')\) is dense. Note that \(\cT(X) \subset \cI(U^x,S^y)\), hence \(\cT(X') = \cT(X) \cap \cI(U^x,S^y) \cap \cU(U^x,S^y)\). Again by the lemma, \(\cI(U^x,S^y) \cap \cU(U^x,S^y)\) is open. Therefore \(\cT(X')\) is open. \qed
\end{prf}

\section{Intersections and orientations}\label{app_intersect_orient}

Here we bring together the necessary preliminaries from intersection theory, in particular, we treat orientations.

Let \(V\) be an \(n\)-dimensional real vector space. An \emph{orienting multivector} of \(V\) is an element \(o \in \Lambda^nV - \{0\}\). If \(S \subset V\) is a subspace, then a \emph{coorienting multivector} (or \emph{form}) of \(S\) is an orienting multivector of the annihilator \(S^\perp \subset V^*\), or, which is the same, an element \(c \in \Lambda^{\codim S}(V/S)^* - \{0\}\). Indeed, the projection \(p \fc V \to V/S\) induces an injective map \(p^* \fc (V/S)^* \to V^*\) whose image is precisely \(S^\perp\), thus \(p \fc (V/S)^* \simeq S^\perp\). The common dimension of these spaces is \(\dim S^\perp = \codim S\), hence \(\Lambda^{\dim S^\perp}S^\perp \simeq \Lambda^{\codim S}(V/S)^*\) canonically. Two (co)orienting multivectors are equivalent if they differ by a positive scalar multiple. An equivalence class of (co)orienting multivectors is a (co)orientation. Note that if \(n=0\), then \(\Lambda^nV = \Lambda^0V = \R\) canonically, and we can speak of the sign of the orientation.

A transverse intersection of an oriented subspace and a cooriented subspace is again an \emph{oriented} subspace. To see this, we need some linear algebra.

For any real vector space \(W\) and integers \(i,j \geq 0\) there are the canonical maps \(\wedge \fc \Lambda^iW \otimes \Lambda^j W \to \Lambda^{i+j}W\) and \(\det \fc \Lambda^iW^* \otimes \Lambda^iW \to \R\). If \(W\) is finite-dimensional, then the determinant defines a non-degenerate pairing between \(\Lambda^iW^*\) and \(\Lambda^iW\), which gives rise to a canonical isomorphism \(\Lambda^iW^* = (\Lambda^iW)^*\).

Let again \(V\) be an \(n\)-dimensional real vector space and let \(0 \leq k \leq l \leq n\) be two integers. Consider the map
\[\Lambda^lV \otimes \Lambda^{k}V^* \otimes \Lambda^{l-k}V^* \xrightarrow{\det \circ (\id_{\Lambda^lV} \otimes \wedge)} \R\,,\]
which can be viewed as a map
\[\frown \fc \Lambda^lV \otimes \Lambda^{k}V^* \to (\Lambda^{l-k}V^*)^* = \Lambda^{l-k}V \,.\]

Now let \(U \subset V\) be an oriented subspace of dimension \(l\) with \(o \in \Lambda^lU\) as an orienting multivector let \(S\) be a cooriented subspace of codimension \(k\) with \(c \in \Lambda^k(V/S)^*\) as a coorienting multivector. Suppose that they intersect transversely. Then \(k \leq l \leq n\). Let \(p \fc U \to V/S\) be the natural projection. Then the element \(o \frown p^*c \in \Lambda^{l-k}U\) can be seen to belong to \(\Lambda^{l-k}(U \cap S) -\{0\}\), thus defining an orienting multivector of \(U \cap S\), hence an orientation.

Note that if \(U \subset V\) is oriented and \(V = U \oplus S\), then on \(S\) there is a natural coorientation. Indeed, let \(o \in \Lambda^kU\) be an orienting multivector. Since \(\Lambda^kU\) is one-dimensional, there is a unique \(\eta \in (\Lambda^kU)^* = \Lambda^kU^*\) with \(\eta(o) = 1\). The projection \(p \fc U \to V/S\) is an isomorphism, hence also \(p^* \fc \Lambda^k(V/S)^* \to \Lambda^kU^*\). The desired coorientation on \(S\) is the the class of the multivector \((p^*)^{-1}(\eta)\). This coorientation on \(S\) can also be defined as the unique coorientation such that the oriented intersection \(U \cap S = \{0\}\) has the positive orientation.

This can be generalized to a transverse intersection of an oriented subspace with a finite number of cooriented ones: if \(U \subset V\) is oriented, \(S_1,\dots,\linebreak S_k \subset V\) are cooriented, then define by induction
\[\textstyle U \cap \bigcap_{i=1}^k S_i := \big(U \cap \bigcap_{i=1}^{k-1}S_i\big)\cap S_k\,.\]

Now if \(M\) is a manifold and \(N,P\) are submanifolds, such that \(N\) is oriented by a continuous section \(o\) of the bundle \(\Lambda^lTN - \{0\}\) (\(l = \dim N\)), and \(P\) is cooriented by \(c \in \Gamma(\Lambda^k (TM/TP)^* - \{0\})\) (\(k = \codim_M P\)), then the intersection \(N \cap P\) acquires a natural orientation, according to the above construction, provided that \(N\) and \(P\) intersect transversely. Again, the transverse intersection of a single oriented submanifold and a finite number of cooriented ones can be given an orientation. Also remark that a simply connected manifold can be (uniquely) oriented by orienting the tangent space to it at some point.

Consider now a Morse-Smale triple \((M,F,\rho)\). By Remark \ref{CW_char_maps} the unstable manifolds are oriented. For \(x \in \Crit(F)\) we have \(T_xM = T_xU^x \oplus T_xS^x\), hence, by the above remark, all the stable manifolds acquire natural coorientations. These coorientations serve to define the intersection numbers introduced in Section \ref{Definitions}, as follows. Let $x\in \Crit_k(F)$ and let $\phi \fc \R^k \to M$ be a smooth map, which is transverse to $S^x$ and satisfies $(\overline{S^x}-S^x) \cap \phi(\R^k) = (\overline{\phi(\R^k)}-\phi(\R^k)) \cap S^x = \varnothing$. In this case the set $\phi^{-1}(S^x)$ is finite. For a point $z \in \phi^{-1}(S^x)$ we define its sign as follows. The differential $d_z\phi$ maps $T_z\R^k = \R^k$ isomorphically onto its image; therefore inside $T_{\phi(z)}M$ we have two subspaces which intersect transversely and which are of complementary dimensions: $d_z\phi(\R^k)$ and $T_{\phi(z)}S^x$. The first space is oriented while the second one is cooriented, thus their intersection, being an oriented subspace of dimension zero, has a number $\ve_z$ associated to it. Now put
$$\#\phi \cap S^x:=\sum_{z\in\phi^{-1}(S^x)}\ve_z\,.$$

Additionally, we would like to remark that in existing literature on Morse homology the spaces of gradient trajectories are oriented in a certain way, see for example \cite{Morse_homology_orient_Hutchings}. Using coorientations, in our opinion, allows for a more canonical (still, equivalent) definition of orientations: for \(x,y \in \Crit(F)\) the manifold \(\cM(x,y) = U^x \cap S^y\) of (parametrized) trajectories between \(x\) and \(y\) carries an orientation as the intersection of an oriented and a cooriented manifolds. The space \(\widetilde \cM(x,y) = \cM(x,y)/\R\) of unparametrized trajectories is oriented because each one of its tangent spaces is a quotient of an oriented space by an oriented subspace.

\section{CW structure from Morse-Smale triples}\label{app_CW_struct}

The following theorem was proven by Paul Biran in \cite{Biran_lagrangian_barriers_symp_embeddings}, based on the results of Laudenbach and Hutchings from \cite{Laudenbach_Thom-Smale_complex} and \cite{Hutchings_Reidemeister_torsion_gend_Morse_thry}. We bring here a more detailed version of his proof, since our results rely on it, and also for the sake of completeness.

\begin{thm}\label{thm_CW_struct} Let \((M,F,\rho)\) be a Morse-Smale triple. Then the unstable manifolds are the open cells of a CW structure on \(M\).
\end{thm}

For the proof we need the following definitions and results.

\begin{defin} An endomorphism of a finite-dimensional real vector space is called hyperbolic if it only has eigenvalues of absolute value different from \(1\).
\end{defin}

\begin{defin}
Let \(X\) be a vector field on \(M\) and let \(\phi_t\) be the flow it generates. A critical point \(x\) of \(X\) is called hyperbolic if \(d_x\phi_t \in \End(T_xM)\) is for all \(t \neq 0\). A closed orbit \(\gamma\) of \(\phi_t\) of period \(T > 0\) is called hyperbolic if the differential of the local Poincar\'e map is.\footnote{The local Poincar\'e map at \(x = \gamma(0)\) is defined as follows: choose a hypersurface \(H\) passing through \(x\) and intersecting \(\gamma\) transversely; now send a point \(y \in H\) sufficiently close to \(x\) to \(\phi_{t(y)}(y)\), where \(t(y) = \inf\{t > 0\,|\phi_t(y) \in H\,\}\).}
\end{defin}

\begin{defin}Let \(X\) be a vector field generating a flow \(\phi_t\). If \(x\) is a hyperbolic critical point, the stable and unstable manifolds through \(x\) are
\[U^x = \{y \in M\,|\, \lim_{t\to-\infty}\phi_ty = x\} \text{ and }S^x = \{y \in M\,|\, \lim_{t\to+\infty}\phi_ty = x\}\,.\]
If \(\gamma\) is a hyperbolic periodic orbit of \(X\) of positive period, similarly we define its stable and unstable manifolds as
\[U^\gamma = \{y \in M\,|\, \exists \tau \in \R: \lim_{t\to-\infty}d(\phi_ty,\gamma(t+\tau)) = 0\}\]
and
\[S^\gamma = \{y \in M\,|\, \exists \tau \in \R: \lim_{t\to+\infty}d(\phi_ty,\gamma(t+\tau)) = 0\}\,,\]
where \(d\) is a metric inducing the topology of \(M\).
\end{defin}

It is well-known that the stable and unstable manifolds of hyperbolic critical points and closed orbits are smooth submanifolds of \(M\).

\begin{defin}A vector field is called Morse-Smale if it only has hyperbolic critical points and closed orbits and moreover, any stable manifold meets any unstable manifold transversely. A Morse-Smale field is called gradient-like if it has no closed orbits of positive period. A gradient-like Morse-Smale field is said to be of standard type if around each critical point there is a coordinate system \((x_1,\dots,x_n)\) such that the field is given by
\[X(x_1,\dots,x_n) = (-x_1, \dots, - x_k, x_{k+1},\dots,x_n)\,.\]
\end{defin}

\begin{defin} Let \((M,F,\rho)\) be a Morse-Smale triple. It is said to be of standard type if the vector field \(-\nabla_\rho f\) is of standard type.
\end{defin}

The following result is proved by Laudenbach \cite{Laudenbach_Thom-Smale_complex} and Hutchings \cite{Hutchings_Reidemeister_torsion_gend_Morse_thry}.
\begin{thm}\label{thm_CW_struct_vf_std_form} Let \((M,F,\rho)\) be a Morse-Smale triple of standard type. Then the unstable manifolds are the open cells of a CW structure on \(M\). \qed
\end{thm}

\begin{thm}[Smale \cite{Smale_grad_dyn_systems}]\label{thm_Smale_grad_like_MS_vf}A gradient-like Morse-smale vector field of standard type is the negative gradient field of a Morse function with respect to some Riemannian metric on \(M\). \qed
\end{thm}

\begin{thm}[Franks \cite{Franks_MS_flows_homotopy_thry}]If \(X_0\) is a gradient-like Morse-Smale vector field, then there exists a continuous path \([0,1]\ni t \mapsto X_t\) in the space of smooth vector fields on \(M\) starting at \(X_0\) such that every \(X_t\) is Morse-Smale and \(X_1\) is a gradient-like Morse-Smale vector field of standard type whose critical points coincide with those of \(X_0\). \qed
\end{thm}

\begin{defin}Two vector fields \(X,Y\) on \(M\) are topologically conjugate if there is a homeomorphism \(\phi\) of \(M\) which maps any orbit of the flow of \(X\) onto an orbit of the flow of \(Y\) and preserves their orientation.
\end{defin}

If \(X\) and \(Y\) are topologically conjugate gradient-like Morse-Smale vector fields, then the above homeomorphism maps critical points of \(X\) onto those of \(Y\), preserves their indices, and also the stable and unstable manifolds.

\begin{defin}A vector field is called structurally stable if any sufficiently \(C^\infty\)-small perturbation of it is topologically conjugate to it.
\end{defin}

The following result was probably first proved by Smale, but in his original paper we could not find quite the formulation we need, so instead we refer the reader to the paper by Robinson \cite{Robinson_struct_stab_vect_fields}.
\begin{thm}Any Morse-Smale vector field is structurally stable. \qed
\end{thm}

\begin{cor}If \(\{X_t\}_{t\in[0,1]}\) is a continuous path of Morse-Smale vector fields then \(X_0\) and \(X_1\) are topologically conjugate. \qed
\end{cor}

\begin{cor}Any gradient-like Morse-Smale vector field is topologically conjugate to one of standard form.\qed
\end{cor}

\begin{prf}[of Theorem \ref{thm_CW_struct}] Let \(X = -\nabla_\rho F\). By the last corollary, \(X\) is topologically conjugate by a homeomorphism \(\phi\) to a gradient-like Morse-Smale vector field \(Y\) of standard form, which, by Theorem \ref{thm_Smale_grad_like_MS_vf}, is the negative gradient vector field of a Morse function with respect to some Riemannian metric. Theorem \ref{thm_CW_struct_vf_std_form} says that the unstable manifolds of \(Y\) are the open cells of a CW structure on \(M\). This means that the unstable manifolds of \(X\) are the open cells of another CW structure on \(M\), which is obtained from that corresponding to \(Y\) by applying \(\phi^{-1}\). \qed
\end{prf}

\section{Morse complex and the chain complex of the CW structure associated to a Morse-Smale triple}\label{app_Morse_complex}

The goal of this section is to prove that the Morse complex of a Morse-Smale triple is canonically isomorphic to the chain complex of the CW structure associated to the triple by the previous section. In order to do this we need to show that if \(x\) is a critical point of index \(k\) and \(y\) one of index \(k-1\), then the algebraic number of connecting trajectories from \(x\) to \(y\) coincides with the degree of map \(S^{k-1} \to M^{k-1}\to M^{k-1}/(M^{k-1}-U^y)\to (\ol {U^y}\cup M^{k-2})/M^{k-2} \simeq S^{k-1}\), where the first arrow is the restriction of the attaching map of the cell corresponding to $U^x$ to the boundary $S^{k-1} = \partial D^k$.

More precisely, we shall prove the following theorem:
\begin{thm}Let \((M,F,\rho)\) be a Morse-Smale triple. Let \(\Phi_x \fc D^k \to M\) be the attaching map for a critical point \(x\in\Crit_k(F)\), given by the CW structure from the last section. These maps induce orientations on the unstable manifolds of the system. Then for any critical points \(x,y\) of indices \(k,k-1\) we have \(\#\widetilde \cM(x,y) = \deg (p_y \circ \Phi_x|_{\partial D^k})\), where \(p_y \fc M^{k-1} \to M^{k-1}/(M^{k-1}-U^y)\simeq S^{k-1}\), the latter sphere carrying the orientation induced from that of $U^y$.
\end{thm}

\begin{prf}Coorient the stable manifolds as described in Appendix \ref{app_intersect_orient}.

Choose a coordinate system near \(y\), \(\Psi\fc B^{k-1} \times B^{n-k+1} \to M\), where \(B^i \subset \R^i\) is a small open ball centered at the origin. Assume \(\Psi(0) = y\) and \(\Psi(B^{k-1} \times 0) \subset U^y\), \(\Psi(0 \times B^{n-k+1}) \subset S^y\). Denote \(D_y = \Psi(B^{k-1} \times 0)\), \(T_y = \im\Psi\) and define \(\pi_y \fc T_y \to D_y\) by \(\pi_y = \Psi \circ \text{pr}_1 \circ \Psi^{-1}\), where \(\text{pr}_1 \fc B^{k-1} \times B^{n-k+1} \to B^{k-1} = B^{k-1} \times 0\) is the projection onto the first factor.

Let \(p_y' \fc M^{k-1} \to M^{k-1}/(M^{k-1} - D_y)\) be the quotient map. Recall that the latter space is a sphere, oriented as \(U^y\). Note that \(p_y'\) is the composition of \(p_y\) and the natural quotient map \(M^{k-1}/(M^{k-1}-U^y) \to M^{k-1}/(M^{k-1}-D_y)\), which is a homotopy equivalence. Thus \(\deg (p_y \circ \Phi_x|_{\partial D^k}) = \deg (p'_y \circ \Phi_x|_{\partial D^k})\).

Recall that if \(S\) is a metric sphere, then \(\eta(S) > 0\) is a number such that any two \(\eta(S)\)-close maps from a compact space into \(S\) are homotopic. Let \(S_y = M^{k-1}/(M^{k-1}-D_y)\) and \(\eta = \eta(S_y)\).

Let \(\pi \fc D^k - 0 \to S^{k-1}\) be the radial projection. Via the attaching map \(\Phi\) it induces a deformation retraction \(\widetilde \pi \fc \Phi(D^k-0) \to \Phi(S^{k-1}) \subset M^{k-1}\). For \(\ve > 0\), let \(\delta(\ve)>0\) be such that if \(z\in\Phi(D^k-0) \cap M^{k-1}_{\delta(\ve)}\), then \(d(\widetilde \pi(z),z)<\ve\).

Consider now a closed ball \(D\) in \(U^x\) of small radius around \(x\). It is diffeomorphic to the standard ball \(D^k\). Choose an orientation-preserving diffeomorphism \(D^k \to D\) and denote by \(\nu\) its restriction to \(\partial D^k = S^{k-1}\). Then \(\nu\) is transverse to all gradient trajectories. It is easy to see that \(\#\widetilde \cM(x,y)\) equals the intersection number of \(f^t\circ\nu\) with \(S^y\) for any \(t\).

There exists \(t_0\) such that for all \(t\geq t_0\) the map \(\pi_x^t:=\pi \circ (\Phi_x|_{U^x})^{-1} \circ f^t \circ \nu\) has degree \(1\). Hence \(\deg (p_y' \circ \Phi_x|_{\partial D^k}) = \deg (p_y' \circ \Phi_x|_{\partial D^k} \circ \pi^t_x)\) for all \(t \geq t_0\). Note that \(p_y' \circ \Phi_x|_{S^{k-1}}\circ \pi_x^t = p_y' \circ \widetilde \pi \circ f^t \circ \nu \fc S^{k-1} \to S_y\).

For any \(\ve>0\) there is \(T_\ve > 0\) such that for any \(t \geq T_\ve\) we have \(\im(f^t\circ\nu) \subset M^{k-1}_\ve\).

Choose now \(\ve > 0\) such that \(\ve < \delta(\eta/2)\), \(T_y\cup(M^{k-1}-D_y)_\ve \supset M^{k-1}_\ve\) and \(d(\pi_y(z),z)<\eta/2\) for \(z \in T_y \cap M^{k-1}_\ve\). Then there is a continuous map \(\pi_y' \fc M^{k-1}_\ve \to S_y\) defined by \(\pi_y'(z) = p_y'(\pi_y(z))\) if \(z\in T_y \cap M^{k-1}_\ve\) and \(\pi_y'(z)=[M^{k-1} - D_y]\) otherwise.

For \(t \geq T_\ve\) there is a map \(\widetilde \nu_t \fc S^{k-1} \to S_y\), defined by \(\widetilde \nu^t = \pi_y' \circ f^t\circ\nu\). The distance between \(\widetilde\nu^t\) and \(p_y' \circ \widetilde \pi \circ f^t \circ \nu\) is less than \(\eta\) and so they are homotopic. Consequently \(\deg(p_y\circ\Phi_x|_{S^{k-1}}) = \deg(p_y' \circ \widetilde \pi \circ f^t \circ \nu) = \deg(\widetilde\nu^t)\). This last degree can be computed at the regular value \(y\) of \(\widetilde\nu^t\) and it is easily seen to be equal to the intersection number of \(\nu\) with \(S^y\). \qed
\end{prf}

\end{document}